\documentclass[10pt]{article}
\usepackage[active]{srcltx}
\usepackage{multirow}
\usepackage{latexsym}
\usepackage{amsmath,amssymb,amsthm,amsfonts}
\usepackage{graphicx}
\usepackage{wrapfig}
\usepackage{hyperref}
\usepackage{mdframed}

\usepackage{mathtools,amscd}

\newcommand{\btd}{\bigtriangledown}
\newcommand{\btu}{\bigtriangleup}
\newcommand{\ba}{\begin{array}}
\newcommand{\ea}{\end{array}}
\newcommand{\half}{\mbox{\scriptsize $\frac{1}{2}$} }

\newcommand{\dst}{\displaystyle}
\newcommand{\refe}[1]{(\ref{#1})}
\newcommand{\abs}[1]{\left\vert#1\right\vert}

\newcommand{\supp}{\mbox{\rm\,supp\,}}

\newtheorem{Theorem}{Theorem}[section]

\newtheorem{Lemma}[Theorem]{Lemma}
\newtheorem{Proposition}[Theorem]{Proposition}
{ \theoremstyle{definition}
\newtheorem{Definition}[Theorem]{Definition}

\newtheorem{Remark}[Theorem]{Remark} }
\numberwithin{equation}{section}

\begin{document}

\title{Multiple $q$-Kravchuk polynomials}
\author{J. Arves\'u\thanks{The research of J. Arves\'u was funded by Agencia Estatal de Investigaci\'on of Spain, 
grant number PGC-2018-096504-B-C33 and  Comunidad Aut\'o\-no\-ma de Madrid, grants CC-G07-UC3M/ESP-3339 and CC-G08-UC3M/ESP-4516.}\, and A. M. Ram\'{\i}rez-Aberasturis
 \medskip \\
Department of Mathematics, Universidad Carlos III de Madrid,\\
Avda. de la Universidad, 30, 28911, Legan\'es, Madrid, Spain
}

\date{September 24, 2020}

\maketitle

\begin{abstract} We study a family of type II multiple orthogonal polynomials. We consider orthogonality conditions 
with respect to a vector measure, in which each component is a $q$-analogue of the binomial distribution. The lowering and raising operators as well as the Rodrigues formula for these polynomials are obtained. The difference equation of order $r+1$ is studied. The connection via limit relation between four types of Kravchuk polynomials is discussed.
\end{abstract}



\section{Introduction}

The importance of the eigenfunctions of the generalized equation of hypergeometric type \cite{nsu} is a well-known fact. Recall that among several special eigenfunctions of this equation, the Bessel functions and the classical orthogonal polynomials are widely used in different fields (see for instance the solutions of the Schr\"{o}dinger and Dirac equations \cite{nsu}). Some multiple orthogonal polynomials (enlarging now the set of special functions) have also the property of being eigenfunctions of a differential or difference equation of higher order. In \cite{arvesu-esposito,arvesu-andys1,Coussement-Assche,Kalyagin,lee,Assche_diff-eq} the authors deal with this type of question for differential/difference equations (see also \cite{miki-tsujimoto-vinet-zhedanov,miki-vinet-zhedanov,nda_van_assche} for some developed applications in physics involving these new special functions).

Multiple orthogonal polynomials are related to the simultaneous rational approximation of a set of analytic functions \cite{angelesco,Aptekarev1,Gonchar,Kalyagin1,Kalyagin2,Mahler,kn_Nikishin,Sorokin,Sorokin3}. In this context, several orthogonality conditions appear with respect to a system of measures. Throughout this paper we will use letter $r$ to denote the dimension of both the vector analytic function and the vector measure $\vec{\mu}$, where the components are positive Borel measures supported on a subset of $\mathbb{R}$, with finite moments. In \cite{arvesu_vanAssche} for some discrete vector measures the corresponding multiple orthogonal polynomials were studied. Among several polynomial families, the authors considered the so-called multiple Kravchuk polynomials, whose extension on a non-uniform lattice will be studied in this paper. In particular, the $(r+1)$-order difference equation that has this new extension of multiple Kravchuk polynomials as eigenfunctions will be our main goal.

The content of this paper begins with a preliminary discussion in Section \ref{prelim}, in which the main notations and background materials are addressed. The emphasis will be placed on the difference equations satisfied by the $q$-Kravchuk polynomials and multiple Kravchuk polynomials. These equations are particular cases of a more general difference equation studied in Section \ref{q-Kravchuk-multiple}. This Section \ref{q-Kravchuk-multiple} will deal with a $q$-analogue of multiple Kravchuk polynomials and their algebraic properties, namely the raising operators and Rodrigues-type formula, which gives an explicit expression for the polynomials. Subsection \ref{diff-q-multiple-Krav} contains a detailed study of the $(r+1)$ order difference equation. Finally, in Section \ref{limit-relations-Krav} a connection between four different extensions of Kravchuk polynomials is presented.

\section{Preliminary material}\label{prelim}

Let $x(s)$ be a non-uniform lattice of the discrete variable $s$, $a\leq s\leq b-1$, $0\leq a <+\infty$, $b\in\mathbb{R}^{+}\cup\{+\infty\}$. By $\bigtriangleup y(s)=y(s+1)-y(s)$ and $\bigtriangledown y(s)=\bigtriangleup y(s-1)$ we denote the forward and backward difference operators, respectively.

The second order difference equation 
\begin{equation}
\ba{c} \displaystyle
\sigma(s) \frac{\btu}{\btu x(s-\half)}  \frac{\btd y(s)}{\btd x(s)}
+ \tau(s) \frac{\btu y(s)}{\btu x(s)} + \lambda y(s) =0, \\
y(s)=y(x(s)),\quad \sigma(s)=a_2(x(s)) - \half a_1(x(s)) \btu x(s-\half),
 \quad
\tau(s)=a_1(x(s)),
 \ea
\label{eqdif}
\end{equation}
where $\deg a_{2}\leq
2$, $\deg a_{1}=1$, and $\lambda\in\mathbb{R}$, is a discrete analogue of the hypergeometric equation \cite{nsu}
\begin{equation}
\dst \sigma(x)y''(x)+\tau(x)y'(x) +\lambda
y(x)=0,\quad \deg \sigma\leq
2,\quad \deg \tau=1.
\label{class-hyp-eq}
\end{equation}

From the self-adjoint form of equation \refe{eqdif} follows that its polynomial solutions (also called orthogonal polynomials of a discrete variable) verify the orthogonality relation
\begin{equation}
\displaystyle \sum_{s = a }^{b-1} P_n(x(s))
P_m(x(s))\omega(s) \bigtriangleup x(s-\half) =
\delta_{n,m}||P_n||^2, \label{norm}
\end{equation}
provided that the orthogonalizing weight $\omega(s)$ solves the Pearson-type difference equation (see \cite[pp. 70-72]{nsu}) and the condition 
\begin{equation*}
\dst\left.\sigma(s) \omega(s) x(s-\half)^{k} \right|_{s=a,b} = 0,\quad k= 0,1,\dots,
\end{equation*}
is fulfilled. See \cite{arvesu-qalgebras,cse,floreanini-letourneux-vinet,nsu} for some applications of the aforementioned polynomial solutions of equation \refe{eqdif}.

In Subsection \ref{diff-q-multiple-Krav} we study the reciprocal situation, i.e. to derive a difference equation from the orthogonality relations for some multiple orthogonal polynomials. 

From the point of view of rational approximation, the orthogonal polynomials derived from relation \refe{norm} form the denominator of the Pad\'e approximants $Q_{n}(z)/P_{n}(z)$ to the Cauchy transform of the involved orthogonality measure \cite{Nikishin}. The situation is similar if we have several Cauchy transforms
\begin{equation}
\dst \hat{\mu}_{i}(z)=\int_{{\Omega}_{i}}\frac{d\mu_{i}(x)}{z-x},
\quad z\notin\Omega_{i}\quad i=1,\dots,r,
\label{hat_functions}
\end{equation}
where $\Omega_{i}$ is the smallest interval that contains $\supp(\mu_{i})$ for each vector component of 
$\vec{\mu}=(\mu_{1},\ldots,\mu_{r})$. These components are positive Borel measures with finite moments. They could be continuous or discrete
\begin{gather}
\mu_{i}=\sum\limits_{k=0}^{N_{i}}\omega_{i,k}\delta_{x_{i,k}},
\qquad
\omega_{i,k}>0,
\qquad
x_{i,k}\in \mathbb{R},
\qquad
N_{i}\in \mathbb{N\cup}\{+\infty \},
\qquad
i=1,2,\ldots,r,
\label{dismeasure}
\end{gather}
where $\delta_{x_{i,k}}$ denotes the Dirac delta function and $x_{i_{1},k}\neq x_{i_{2},k}$, $k=0,\ldots,N_{i}$,
whenever $i_{1}\neq i_{2}$. Indeed, the functions \refe{hat_functions} can be simultaneously approximated by rational functions with prescribed order near infinity \cite{Nikishin}
\begin{equation*}
\ba{c}
\dst P_{\vec{n}}(z)\hat{\mu}_{i}(z)-Q_{\vec{n},i}(z)=\frac{\zeta_{i}}{z^{n_{i}+1}}+\cdots=\mathcal{O}(z^{-n_{i}-1}),\quad
i=1,\dots,r,
\ea
\label{f1}
\end{equation*}
where a multi-index $\vec{n}=(n_{1},
n_{2},\dots,n_{r})$ of nonnegative integers is introduced, and a polynomial $P_{\vec{n}}(z)\not\equiv0$ of
degree at most $|\vec{n}|=n_1+\cdots+n_r$ must be found, if any.

Observe that $P_{\vec{n}}(z)$ is a common denominator of the simultaneous rational approximants $Q_{\vec{n},i}(z)/P_{\vec{n}}(z)$, to $\hat{\mu}_{i}(z)$, $i=1,2,\dots,r$. The $P_{\vec{n}}(z)$ is the so-called type II
multiple orthogonal polynomial of degree $\leq|\vec{n}|$ defined by the orthogonality conditions 
\begin{equation}
\int_{\Omega_i}P_{\vec{n}}(x)x^k\ d\mu_i(x) =  0,\quad k=0,1,\ldots,n_i-1,\quad i=1,\dots,r.
\label{stelseltype2}
\end{equation}

The conditions \refe{stelseltype2} give a linear system of $|\vec{n}|$ homogeneous equations for the
$|\vec{n}|+1$ unknown coefficients of $P_{\vec{n}}(z)$. In this paper we will deal with a unique monic polynomial solution of \refe{stelseltype2} with $\deg P_{\vec{n}}(x)=|\vec{n}|$ for every $\vec{n}$.
This situation occurs when the above system of measures forms an $AT$ system~\cite{Nikishin}, that is, every multi-index is normal. We will focus only on such a system of discrete measures, for which $\Omega_{i}=\Omega\subset\mathbb{R}^+$, $i=1,2,\ldots,r$.

In the sequel we will represent any discrete polynomial $P_n(x(s))$ by $P_n(s)$ and consider the non-uniform lattice $x(s)=\frac{q^s-1}{q-1}$, $|q|\not=1$ as well as the following definition \cite{arvesu-qHahn}.

\begin{Definition}
\label{defqpo}
A polynomial $P_{\vec{n}}(s)$ of degree $\vert \vec{n}\vert$ on the lattice $x(s)$ is said to be multiple $q$-orthogonal polynomial of a multi-index $\vec{n} \in \mathbb{N}^{r}$
with respect to positive discrete measures $\mu_{1},\mu_{2},\ldots,\mu_{r}$ (with finite moments), $\supp(\mu_{i})\subset \Omega_{i}\subset \mathbb{R}$, $i=1,2,\ldots,r$, if the following condition holds:
\begin{gather*}
\sum\limits_{s=0}^{N_{i}}P_{\vec{n}}(s)  [s]_{q}^{(k)}d\mu_{i}=0,
\qquad
k=0,\ldots,n_{i}-1,
\qquad
N_{i}\in\mathbb{N\cup}\{+\infty\},
\notag
\end{gather*}
\end{Definition}
where 
\begin{gather*}
[s]_{q}^{(k)}=\prod\limits_{j=0}^{k-1}\frac{q^{s-j}-1}{q-1} =x(s) x(s-1) \cdots x(s-k+1)
\quad\!
\text{for}
\quad
k>0,
\quad\!
\text{and}
\quad\!
[s]_{q}^{(0)}=1, 
\end{gather*}
is the $q$-analogue of the Stirling polynomials.

We will refer to a function $f_q(s)$ as a $q$-analogue to a given function $f(s)$ if for any sequence $(q_n)_{n\geq0}$ approaching to $1$, the corresponding sequence 
$\left(f_{q_{n}}(s)\right)_{n\geq0}$ tends to $f(s)$. In particular, the lattice $x(s)=\frac{q^s-1}{q-1}$ ($s=0,1,\dots,N\in\mathbb{N}$) tends to the uniform lattice $s$, as $q$ approaches to $1$. In addition, by a linear transformation we can transform the lattice $s=0,1,\dots,N$ (with step $1$) into a new uniform lattice with step $h\in\mathbb{R}^{+}$. If $h$ approaches to zero, the equation \refe{eqdif} transforms into the hypergeometric equation \refe{class-hyp-eq} \cite{nsu}.

In \cite{alv_arv} a $q$-analogue of the classical Kravchuk polynomials \cite{kravchuk1,kravchuk2} was studied. For such a purpose the authors considered for the equation \refe{eqdif}, some specific polynomials $a_{2}(s)$ and $a_{1}(s)$ (see \refe{a2_a1_lambda} below) as well as the orthogonality weight function
\begin{equation}
\omega(s)=\left(\frac{p}{1-p}\right)^{s} \frac{ q^{\half s(s-1)} [N]_q!(1-p)^{N} }
{\Gamma_q(N-s+1)\Gamma_q( s+1 ) } ,\quad  0<p<1.
\label{omega_func}
\end{equation}
In \refe{omega_func} the function
 \begin{gather*}
\Gamma_q(s) =
\begin{cases}
f(s;q)=(1-q)^{1-s} \dfrac{\prod\limits_{k\geq0} (1-q^{k+1})}{\prod\limits_{k\geq0} (1-q^{s+k})}, & 0<q<1,
\\
q^{\frac{(s-1)(s-2)}{2}}f\big(s;q^{-1}\big), & q>1,
\end{cases}
\end{gather*}
is a $q$-analogue of the Gamma function \cite{Gasper,nsu} and 
$[N]_q=q^{\frac{1-N}{2}}\dfrac{\Gamma_q(N+1)}{\Gamma_q(N)}$ is a $q$-number.

In addition to the above notations we will use throughout this paper the following difference operators 
\begin{gather}
\Delta  \overset{\rm def}{=}\frac{\bigtriangleup} {\bigtriangleup x(s-1/2)},\qquad
\nabla \overset{\rm def}{=} \frac{\bigtriangledown}{\bigtriangledown x(s+1/2)},\qquad
\nabla^{n_{j}} =\underbrace{\nabla\cdots\nabla}_{n_{j}~\text{times}},\quad n_{j}\in\mathbb{N}.
\label{lower}
\end{gather}
When convenient, the representation  
$\bigtriangledown x_{1}(s)\overset{\rm def}{=}\bigtriangledown x(s+1/2)= \bigtriangleup x(s-1/2)=q^{s-1/2}$ will be used. 

Recall that 
\begin{align}
\bigtriangledown^{m}\left(f(s)g(s)\right)&=\sum_{k=0}^{m}\binom{m}{k}\left(\bigtriangledown^{k}f(s)\right)\left(
\bigtriangledown^{m-k}g(s-k)\right),\quad m\in\mathbb{N},\notag\\
\bigtriangledown^{m}f(s)&=\sum_{k=0}^{m}(-1)^{k}\binom{m}{k}f(s-k).\label{Leibniz2}
\end{align}

\subsection{Multiple Kravchuk polynomials}

Let $N\in\mathbb{N}$ be a parameter and $\vec{n}\in \mathbb{N}^{r}$ be a multi-index with $\left\vert \vec{n}\right\vert \leq N$, and $\vec{p}=\left( p_{1},\ldots ,p_{r}\right)$, where $0<p_{i}<1$, $i=1,2,\ldots ,r$,\ and with all the $p_{i}$ different. Multiple Kravchuk polynomials $K_{\vec{n}}^{\vec{p},N}(x)$ are the unique monic polynomials  of degree $\left\vert \vec{n}\right\vert $ that satisfy the orthogonality conditions
\begin{equation*}
\sum\limits_{x=0}^{N}K_{\vec{n}}^{\vec{p},N}( x) (
-x) _{j}\upsilon ^{p_{i},N}( x) =0,
\qquad
j=0,\ldots,n_{i}-1,
\qquad
i=1,\ldots,r, 
\end{equation*}%
where $(-x)_{j}$ denotes the Pochhammer symbol \cite{Gasper,nsu} and
\begin{equation*}
\upsilon ^{p_{i},N}( x) =
\begin{cases}
\displaystyle\frac{N!p_{i}^{x}\left( 1-p_{i}\right) ^{N-x}}{\Gamma (
x+1) \Gamma \left( N-x+1\right) }, & \text{ if }x\in \mathbb{R}\backslash \left( \mathbb{Z}^{-}\cup \left\{ N+1,N+2,\ldots \right\} \right)
,\cr\cr0, & \text{otherwise.}
\end{cases}
\end{equation*}

In \cite{arvesu_vanAssche} the normality of the the multi-index $\vec{n}\in \mathbb{N}^{r}$ was addressed and the following raising operators were found 
\begin{equation}
\mathcal{L}^{p_{i},N}\left[ K_{\vec{n}}^{\vec{p},N}( x) \right] =-K_{\vec{n}+\vec{e_{i}}}^{\vec{p},N+1}( x),\quad i=1,\ldots ,r, \label{raisingxKra}
\end{equation}
where 
\begin{equation*}
\mathcal{L}^{p_{i},N} \overset{\rm def}{=}\frac{p_{i}\left( 1-p_{i}\right) \left( N+1\right) }{\upsilon
^{p_{i},N+1}( x) }\bigtriangledown \upsilon
^{p_{i},N}( x).
\end{equation*} 
As a~consequence of~\eqref{raisingxKra} the Rodrigues-type formula can be obtained
\begin{equation}
K_{\vec{n}}^{\vec{p},N}(x) =\left( -N\right) _{\left\vert \vec{n}\right\vert }\left( \prod\limits_{i=1}^{r}p_{i}^{n_{i}}\right) \frac{\Gamma(x+1)\Gamma (N-x+1)}{N!} \mathcal{K}_{\vec{n}}^{\vec{p}}\frac{\left( N-\left\vert \vec{n}\right\vert \right) !}{\Gamma (x+1)\Gamma (N-\left\vert \vec{n}\right\vert
-x+1)},  \label{RFK}
\end{equation}
where
\begin{equation}
\mathcal{K}_{\vec{n}}^{\vec{p}}=\prod_{i=1}^{r}\left( \frac{1-p_{i}}{p_{i}}\right) ^{x}\bigtriangledown ^{n_{i}}\left( \frac{p_{i}}{1-p_{i}}\right)
^{x}.\label{operator_Krav_Rodrig}
\end{equation}
In \cite{lee} the author found the high-order linear difference equation 
\begin{gather}
\prod\limits_{i=1}^{r}\mathcal{L}^{p_{i},N+r-i-1}\big[\bigtriangleup K_{\vec{n}}^{\vec{p},N}(x) \big]=-\sum\limits_{i=1}^{r}n_{i}\prod\limits_{\substack{j=1\\
j\neq i}}^{r} \mathcal{L}^{p_{j},N+r-j-1}\big[K_{\vec{n}}^{\vec{p},N}(x)\big].
\label{opdi-1Kra}
\end{gather}
Moreover, the recurrence relation was found in \cite{arvesu_vanAssche,Hane}
\begin{multline}
xK_{\vec{n}}^{\vec{p},N}( x)=K_{\vec{n}+\vec{e}_{k}}^{\vec{p}
,N}( x) +\left[ (N-\left\vert \vec{n}\right\vert
)p_{k}+\sum\limits_{i=1}^{r}n_{i}\left( 1-p_{i}\right) \right] K_{\vec{n}}^{
\vec{p},N}( x)\\
+\sum\limits_{i=1}^{r}n_{i}p_{i}\left( p_{i}-1\right) \left( \left\vert \vec{n}\right\vert -N-1\right) K_{\vec{n}-\vec{e}_{i}}^{\vec{p},N}(x),
 \label{RRKra}
\end{multline}
where $\vec{e}_{j}$ denotes the standard $r$-dimensional unit vector with 
the $j$-th entry equals $1$ and $0$ otherwise, $1\leq j\leq r$.

Observe that the multiple Kravchuk polynomials $K_{\vec{n}}^{\vec{p},N}(x)$ are common eigenfunctions of the $(r+1)$-order linear difference equations \eqref{opdi-1Kra} and ~\eqref{RRKra}.

\section{Multiple $q$-Kravchuk polynomials}\label{q-Kravchuk-multiple}
Consider a multi-index $\vec{n}\in \mathbb{N}^{r}$ and the following~$r$ positive discrete measures on a subset of $\mathbb{R}^{+}$,
\begin{gather}
\mu_{i}=\sum\limits_{s=0}^{N}\omega_{i}(k)\delta(k-s),
\qquad
\omega_{i}>0,
\qquad
i=1,2,\ldots,r,
\label{MeasuK2}
\end{gather}
where 

\begin{equation*}
\omega_{i}(s)=\upsilon_{q}^{p_{i}, \beta_{i}, N}(s) \bigtriangleup x(s-1/2),\quad
\upsilon_{q}^{p_{i}, \beta_{i}, N}( s) =
\begin{cases}
\displaystyle\frac{q^{\binom{s}{2}}\left[N\right]_{q}!p_{i}^{s}\beta_{i}^{N-s}}{\Gamma_{q} (
s+1) \Gamma_{q} \left( N-s+1\right) }, & \text{ if }s\in \Theta 
,\cr\cr0, & \text{otherwise,}
\end{cases}
\end{equation*}
with $\Theta = \mathbb{R}\backslash \left( \mathbb{Z}^{-}\cup \left\{ N+1,N+2,\ldots \right\} \right)$, $\left\vert \vec{n}\right\vert \leq N$, $\beta_{i}=\left( 1-p_{i}\right)$, $0<p_{i}<1$, $i=1,2,\ldots ,r$, and with all the $p_{i}$ different. As a consequence of \cite[Lemma 5.1]{arvesu-andys2} with $\alpha_{i}=\frac{p_i}{1-p_i}$, the system of measures $\mu_1,\mu_2,\ldots,\mu_r$ given in~\eqref{MeasuK2} forms an AT system on $\Theta$. Next, we will used Definition \ref{defqpo} with respect to the measures~\eqref{MeasuK2}.

\begin{Definition}
A polynomial $K_{q,\vec{n}}^{\vec{p},\vec{\beta},N}(s)$ of degree $\vert\vec{n}\vert$ ($\vec{n}\in \mathbb{N}^{r}$) that verifies the orthogonality conditions 
\begin{gather}
\sum\limits_{s=0}^{N}K_{q,\vec{n}}^{\vec{p},\vec{\beta},N}(s) [s]_{q}^{(k)}\upsilon_{q}^{p_{i}, \beta_{i}, N}(s)
\bigtriangleup x(s-1/2) =0,
\qquad
0\leq k\leq n_{i}-1,
\qquad
i=1,\ldots,r,
\label{NOrthCK2}
\end{gather}
is said to be the $q$-Kravchuk multiple orthogonal polynomial.
\label{defin_qkrav}
\end{Definition}
 
In this paper we will consider monic $q$-Kravchuk multiple orthogonal polynomials. When $r=1$ we recover the monic $q$-Kravchuk polynomials computed in \cite{alv_arv} with respect to the aforementioned weight function \refe{omega_func}.

Observe that $K_{q,\vec{n}}^{\vec{p},\vec{\beta},N}(s)$ has exactly $\vert \vec{n}\vert$ different zeros on 
$\mathbb{R}^{+}$ because of the AT-property (see~\cite[Theorem 2.1, pp.~26--27]{arvesu_vanAssche}).

\begin{Lemma}\label{raising_operator_Krav}
For monic $q$-Kravchuk multiple orthogonal polynomials we have $r$ raising operators
\begin{gather}
\mathcal{D}_{q}^{p_{i}, \beta_{i}, N}K_{q,\vec{n}}^{\vec{p},\vec{\beta},N}(s)
=-q^{1/2}K_{q,\vec{n}+\vec{e}_{i}}^{\vec{p},\vec{\beta}_{i,q^{2}},N+1}(s),
\qquad
i=1,\ldots,r,
\label{ROpqMKrav}
\end{gather}
where $\vec{\beta}_{i,q^{2}}=(\beta_{1},\ldots,q^{2}\beta_{i},\ldots,\beta_{r})$ and
\begin{gather}
\mathcal{D}_{q}^{p_{i}, \beta_{i}, N}\overset{\rm def}{=}
\left(\frac{p_{i}\beta_i q^{|\vec{n}|+2N+1}[N+1]_q}{\left[p_{i}\left( q^{|\vec{n}|-1}-1\right)+1\right]\upsilon_{q}^{p_{i},N+1 ,q^{2}\beta_i}(s)} \nabla \upsilon_{q}^{p_{i}, \beta_{i}, N}(s)\right).
\label{D_raising_Krav}
\end{gather}
\end{Lemma}

Notice that we call $\mathcal{D}_{q}^{p_{i}, \beta_{i}, N}$ a raising operator because the $i$-th component of the multi-index
$\vec{n}$ in \eqref{ROpqMKrav} is increased by $1$.

\begin{proof}
From the equation $[s]_{q}^{(k)}=\left(q^{k-1/2}/[k+1]_{q}^{(1)}\right)\nabla[s+1]_{q}^{(k+1)}$,
we rewrite relation \eqref{NOrthCK2} as follows
\begin{gather*}
\sum\limits_{s=0}^{N}K_{q,\vec{n}}^{\vec{p},\vec{\beta},N}(s) \nabla
[s+1]_{q}^{(k+1)}\upsilon_{q}^{p_{i}, \beta_{i}, N}(s) \bigtriangleup x(s-1/2) =0,
\qquad
0\leq k\leq n_{i}-1,
\qquad
i=1,\ldots,r.
\qquad
\end{gather*}
Using summation by parts together with $\upsilon_{q}^{p_{i}, \beta_{i}, N}(-1)=
\upsilon_{q}^{p_{i}, \beta_{i}, N}(N+1)=0$ yields
\begin{multline*}
\sum\limits_{s=0}^{N}\nabla\left( K_{q,\vec{n}}^{\vec{p},\vec{\beta},N}(s)\upsilon_{q}^{p_{i}, \beta_{i}, N}(s)\right) 
[s]_{q}^{(k+1)} \bigtriangleup x(s-1/2) \\
= -\sum\limits_{s=0}^{N}K_{q,\vec{n}}^{\vec{p},\vec{\beta},N}(s) \upsilon_{q}^{p_{i}, \beta_{i}, N}(s)\Delta
[s]_{q}^{(k+1)} \bigtriangleup x(s-1/2)\\
= -\sum\limits_{s=0}^{N}K_{q,\vec{n}}^{\vec{p},\vec{\beta},N}(s) \upsilon_{q}^{p_{i}, \beta_{i}, N}(s)\nabla
[s+1]_{q}^{(k+1)} \bigtriangleup x(s-1/2).
\end{multline*}
Equivalently,
\begin{equation}
\sum\limits_{s=0}^{N+1}\nabla\left( K_{q,\vec{n}}^{\vec{p},\vec{\beta},N}(s)\upsilon_{q}^{p_{i}, \beta_{i}, N}(s)\right) 
[s]_{q}^{(k+1)} \bigtriangleup x(s-1/2)=0,
\qquad
0\leq k\leq n_{i}-1,
\qquad
i=1,\ldots,r,\label{eq_orthog=0}
\end{equation}
where
\begin{equation*}
\nabla\left(K_{q,\vec{n}}^{\vec{p},\vec{\beta},N}(s)\upsilon_{q}^{p_{i}, \beta_{i}, N}(s)\right)=q^{-s+1/2}\left(\upsilon_{q}^{p_{i}, \beta_{i}, N}(s)K_{q,\vec{n}}^{\vec{p},\vec{\beta},N}(s)
-\upsilon_{q}^{p_{i}, \beta_{i}, N}(s-1)K_{q,\vec{n}}^{\vec{p},\vec{\beta},N}(s-1)\right),
\end{equation*}
and
\begin{align*}
\upsilon_{q}^{p_{i}, \beta_{i}, N}(s)&=q^{-s}\frac{x(N+1)-x(s)}{(1-p_i)\left[N+1\right]_q}
\upsilon_{q}^{p_{i},\beta_{i},N+1}(s),\\
\upsilon_{q}^{p_{i}, \beta_{i}, N}(s-1)&=q^{-s}\frac{qx(s)}{p_i\left[N+1\right]_q}
\upsilon_{q}^{p_{i},\beta_{i},N+1}(s).
\end{align*}
Therefore
\begin{equation*}
\nabla\left(K_{q,\vec{n}}^{\vec{p},\vec{\beta},N}(s)\upsilon_{q}^{p_{i}, \beta_{i}, N}(s)\right)=-q^{1/2}
\frac{p_i\left(q^{|\vec{n}|-1}-1\right)+1}
{p_i\beta_iq^{|\vec{n}|+2N+1}\left[N+1\right]_q}
\upsilon_{q}^{p_{i}, q^{2}\beta_{i}, N+1}(s)
\mathcal{P}_{q,\vec{n}+\vec{e}_i}(s),
\end{equation*}
where $\mathcal{P}_{q,\vec{n}+\vec{e}_i}(s)$ denotes a monic polynomial $x^{|\vec{n}|+1}+ \text{lower degree terms}$. Consequently, from \refe{eq_orthog=0}
\begin{equation*}
\sum\limits_{s=0}^{N+1}\mathcal{P}_{q,\vec{n}+\vec{e}_i}(s)\upsilon_{q}^{p_{i}, q^{2}\beta_{i}, N+1}(s) 
[s]_{q}^{(k+1)} \bigtriangleup x(s-1/2)=0.
\end{equation*}
From the uniqueness of the polynomial system defined by \eqref{NOrthCK2} we have that $\mathcal{P}_{q,\vec{n}+\vec{e}_i}(s)=K_{q,\vec{n}+\vec{e}_{i}}^{\vec{p}, \vec{\beta}_{i,q^{2}}, N+1}(s)$, i.e.
\begin{equation*}
\nabla\left(K_{q,\vec{n}}^{\vec{p},\vec{\beta},N}(s)\upsilon_{q}^{p_{i}, \beta_{i}, N}(s)\right)=-q^{1/2}
\frac{p_i\left(q^{|\vec{n}|-1}-1\right)+1}
{p_i\beta_iq^{|\vec{n}|+2N+1}\left[N+1\right]_q}
\upsilon_{q}^{p_{i}, q^{2}\beta_{i}, N+1}(s)K_{q,\vec{n}+\vec{e}_{i}}^{\vec{p}, \vec{\beta}_{i,q^{2}}, N+1}(s),
\end{equation*}
which is equivalent to \refe{ROpqMKrav}.
\end{proof}

\begin{Proposition}
The following $q$-analogue of the Rodrigues-type formula holds
\begin{equation}
K_{q,\vec{n}}^{\vec{p},\vec{\beta},N}(s)
=\mathcal{G}_{q}^{\vec{n},\vec{p},N}\frac{\Gamma_{q}(N-s+1)\Gamma_{q}(s+1)}{q^{\binom{s}{2}}[N]_q!}
\prod\limits_{i=1}^{r}\left(\frac{\beta_i}{p_{i}}\right)^{s} \nabla^{n_{i}}\left(\frac{p_{i}}{\beta_i}\right)^{s}
\frac{q^{\binom{s}{2}+2n_{i}s}[N-|\vec{n}|]_q!}{\Gamma_{q}(N-\vert \vec{n} \vert-s+1)\Gamma_{q}(s+1)},
\label{RFormulaKrav}
\end{equation}
where

\begin{gather}
\mathcal{G}_{q}^{\vec{n},\vec{p},N}=(-1)^{\vert \vec{n} \vert}[N]_{q}^{\left(\vert \vec{n} \vert\right)}
q^{-\frac{5\vert \vec{n}\vert}{2}}\left(\prod\limits_{i=1}^{r}\frac{p_i^{n_i}}{\prod\limits_{j=1}^{n_i}q^{-j}\left[p_{i} \left( q^{|\vec{n}|-|\vec{n}|_i-j -1}-1\right)+1\right]}\right)
\left(\prod\limits_{i=1}^{r-1}q^{n_i\sum\limits_{j=i+1}^{r}n_{j}}\right),\label{Rodri_coef}
\end{gather}
and $\vert \vec{n}\vert_{i}=n_{1}+\cdots +n_{i-1}$, $\vert \vec{n}\vert_{1}=0$.
\end{Proposition}

\begin{proof} From Lemma \ref{raising_operator_Krav}, by applying $k_i$-times ($i=1,\dotsc,r$) the raising operators~\eqref{ROpqMKrav} in a recursive way we obtain
\begin{multline*}
\prod_{i=1}^{r}\left(\frac{p_{i}q^{-2k_i}}{1-p_i}\right)^{-s} \nabla^{n_{i}}\left(\frac{p_{i}}{1-p_i}\right)^{s}\frac{q^{\binom{s}{2}}[N]_q!}{\Gamma_{q}(N-s+1)\Gamma_{q}(s+1)}K_{q,\vec{n}}^{\vec{p},\vec{\beta},N}(s)=[N+1]_{q}^{\left(\vert \vec{k} \vert\right)}q^{|\vec{k}|/2}\\
\left(\prod\limits_{i=1}^{r}\frac{\prod\limits_{j=1}^{k_i}\left(\alpha_i q^{|\vec{n}|+\beta-j}-1\right)}{p_i^{k_i}(1-p_i)^{k_i}}\right)\prod_{i=1}^{r}q^{-n_i\sum\limits_{j=i}^{r}k_j}\prod_{i=1}^{r-1}q^{-k_i\sum\limits_{j=i+1}^{r}k_j}\\
K_{q,\vec{n}+\vec{k}}^{p_1,\dotsc,p_r,N+|\vec{k}|,q^{2k_1}\beta_1,\dotsc,q^{2k_r}\beta_r}(s)\frac{q^{\binom{s}{2}}[N+|\vec{k}|]_q!}{\Gamma_{q}(N+|\vec{k}|-s+1)\Gamma_{q}(s+1)}.
\end{multline*}
Taking in the above expression $n_1=n_2=\dotsb =n_r=0$ and replacing $\beta_i$ by $\beta_iq^{-2k_i}$, $N$ by $N-|\vec{k}|$,  and $k_i$ by $n_i$, for $i=1,\dotsc,r$, the formula \eqref{RFormulaKrav} holds.

\end{proof}

\subsection{Difference equation for multiple $q$-Kravchuk polynomials}\label{diff-q-multiple-Krav}

In this section we will find a lowering operator for the $q$-Kravchuk multiple orthogonal polynomials.
Then we will combine it with the raising operators~\eqref{ROpqMKrav} to get an $(r+1)$-order difference equation on the lattice $x(s)$. More specifically we will follow the next steps: 1) By using some interpolation conditions we define a finite subspace $\mathbb{V}_{r}$ of polynomials on $x(s)$ of degree at most $|\vec{n}| - 1$. 2) We express the lowering operator on the polynomials as a linear combination of the basis vectors of $\mathbb{V}_{r}$.  3) We combine the lowering and the raising operators to derive the difference equation that has the multiple $q$-Kravchuk polynomials as eigenfunctions. 

Despite the above description, the calculations involved in all steps depend on the weights, i.e., on the specific family of multiple orthogonal polynomials (see \cite{arvesu-qHahn,arvesu-esposito,arvesu-andys} for some $q$-multiple orthogonal polynomials, \cite{lee,Assche_diff-eq} for discrete multiple orthogonal polynomials, and \cite{Coussement-Assche} for classical multiple orthogonal polynomials). 

\begin{Lemma}
Let $\mathbb{V}_{r}$ be the linear subspace of polynomials $Q(s)$ on the lattice $x(s)$ of degree at most $\vert
\vec{n}\vert-1$ defined by the following conditions
\begin{gather*}
\sum\limits_{s=0}^{N}Q(s) [s]_{q}^{(k)}\upsilon_{q}^{p_{j},\beta_j/q^{2},N-1}(s)\bigtriangledown x_{1}(s)=0,
\qquad
0\leq k\leq n_{j}-2
\qquad
\text{and}
\qquad
j=1,\ldots,r.
\end{gather*}
Then, the system $\left\{K_{q,\vec{n}-\vec{e}_{i}}^{\vec{p},\vec{\beta}_{i,1/q^{2}},N-1}(s) \right\}_{i=1}^{r}$, where
$\vec{\beta}_{i,1/q^{2}}=(\beta_{1},\ldots,\beta_{i}/q^{2},\ldots,\beta_{r})$, is a~basis for $\mathbb{V}_{r}$.
\label{LIK2}
\end{Lemma}

\begin{proof}

The polynomials $K_{q,\vec{n}-\vec{e}_{i}}^{\vec{p},\vec{\beta}_{i,1/q^{2}},N-1}(s)$, $i=1,\ldots,r$, 
verify the orthogonality relations
\begin{gather*}
\sum\limits_{s=0}^{N}K_{q,\vec{n}-\vec{e}_{i}}^{\vec{p},\vec{\beta}_{i,1/q^{2}},N-1}(s)
[s]_{q}^{(k)}\upsilon_{q}^{p_{j},\beta_j/q^{2},N-1}(s) \bigtriangledown x_{1}(s) =0,
\qquad
0\leq k\leq n_{j}-2,
\qquad
j=1,\ldots,r,
\end{gather*}
that is, they belong to $\mathbb{V}_{r}$.

Assume that there exist numbers $\lambda_{i}$, $i=1,\ldots,r$, such that
\begin{gather}
\sum\limits_{i=1}^{r}\lambda_{i}K_{q,\vec{n}-\vec{e}_{i}}^{\vec{p},\vec{\beta}_{i,1/q^{2}},N-1}(s) =0,
\qquad
\text{where}
\qquad
\sum\limits_{i=1}^{r}\abs{\lambda_{i}}>0.\label{eq_wt_lambda}
\end{gather}
Multiplying \refe{eq_wt_lambda} by $[s]_{q}^{(n_{k}-1)}\upsilon_{q}^{p_{k},\beta_k,N}(s)\bigtriangledown
x_{1}(s)$ and summing on~$s$ from $0$ to~$N$, one has
\begin{gather*}
\sum\limits_{i=1}^{r}\lambda_{i}\sum\limits_{s=0}^{N}K_{q,\vec{n}-\vec{e}_{i}}^{\vec{p},\vec{\beta}_{i,1/q^{2}},N-1}(s)[s]_{q}^{(n_{k}-1)}\upsilon_{q}^{p_{k},\beta_k,N}(s)\bigtriangledown x_{1}(s) =0.
\end{gather*}
Thus, taking into account that
\begin{gather}
\sum\limits_{s=0}^{N}K_{q,\vec{n}-\vec{e}_{i}}^{\vec{p},\vec{\beta}_{i,1/q^{2}},N-1}(s)[s]_{q}^{(n_{k}-1)}\upsilon_{q}^{p_{k},\beta_k,N}(s) \bigtriangledown x_{1}(s) =c\delta_{i,k},
\label{Orthog_CondK2}
\end{gather}
where $\delta_{i,k}$ is the Kronecker delta symbol and $c$ denotes a non-zero real number, we deduce that $\lambda_{k}=0$ for $k=1,\ldots,r$, which contradicts our assumption in \refe{eq_wt_lambda}. This implies that $\{K_{q,\vec{n}-\vec{e}_{i}}^{\vec{p},\vec{\beta}_{i,1/q^{2}},N-1}(s) \}_{i=1}^{r}$ is linearly independent in
$\mathbb{V}_{r}$.
Furthermore, we know that any polynomial of $\mathbb{V}_{r}$ can be determined with $\vert \vec{n}\vert$
coef\/f\/icients while $(\vert \vec{n}\vert -r)$ linear conditions are imposed on $\mathbb{V}_{r}$, i.e.,
$\dim\mathbb{V}_{r}\leq r$.
Therefore, the system $\{K_{q,\vec{n}-\vec{e}_{i}}^{\vec{p},\vec{\beta}_{i,1/q^{2}},N-1}(s)\}_{i=1}^{r}$ spans $\mathbb{V}_{r}$, which completes the proof.
\end{proof}

Now we will prove that the operator~\eqref{lower} is indeed a~lowering operator for the sequence of $q$-Kravchuk multiple orthogonal polynomials $K_{q,\vec{n}}^{\vec{p},\vec{\beta},N}(s)$. The expression in Lemma \ref{lemma_lower} can be viewed as a lowering operator since the degree is lowered by one.

\begin{Lemma}
The following expansion holds
\begin{gather}
\Delta K_{q,\vec{n}}^{\vec{p},\vec{\beta},N }(s) =\sum\limits_{i=1}^{r}q^{|\vec{n}|-n_{i}+1/2}
\left(\frac{ p_{i}\left( q^{n_{i}}-1\right)+1}{ p_{i}\left( q^{|\vec{n}|}-1\right)+1}\right)
[n_{i}]_{q}^{(1)}K_{q,\vec{n}-\vec{e}_{i}}^{\vec{p},\vec{\beta}_{i,1/q^{2}},N-1}(s).
\label{Rela_qK2}
\end{gather}
\label{lemma_lower}
\end{Lemma}

\begin{proof}
Using summation by parts we have the following relation
\begin{align}
\sum\limits_{s=0}^{N}\Delta K_{q,\vec{n}}^{\vec{p},\vec{\beta},N }(s)[s]_{q}^{(k)}
\upsilon_{q}^{p_{j},\beta_j/q^{2},N-1}(s) \bigtriangledown x_{1}(s)
&=-\sum\limits_{s=0}^{N}K_{q,\vec{n}}^{\vec{p},\vec{\beta},N }(s)\nabla\left([s]_{q}^{(k)}\upsilon_{q}^{p_{j},\beta_j/q^{2},N-1}(s)\right) \bigtriangledown x_{1}(s)
\nonumber
\\
&=-\sum\limits_{s=0}^{N}K_{q,\vec{n}}^{\vec{p},\vec{\beta},N }(s)\varphi_{j,k}(s)\upsilon_{q}^{p_{j},\beta_{j},N}(s) \bigtriangledown x_{1}(s), \label{inte-K12}
\end{align}
where
\begin{gather*}
\varphi_{j,k}(s) =q^{1/2}\left(-\frac{q^{-2(N-1)}x(s)}{(1-p_j)[N]_q}+\frac{q^{-2(N-1)}x(N)}{(1-p_j)[N]_q}\right)[s]_{q}^{(k)}-q^{-1/2}\frac{q^{-2(N-1)}x(s)}{p_{j}[N]_q}[s-1]_{q}^{(k)},
\end{gather*}
is a~polynomial of degree $\leq k+1$ in the variable $x(s)$. Therefore, from the orthogonality conditions~\eqref{NOrthCK2}
\begin{gather*}
\sum\limits_{s=0}^{N}\Delta K_{q,\vec{n}}^{\vec{p},\vec{\beta},N }(s) [s]_{q}^{(k)}\upsilon_{q}^{p_{j},\beta_j/q^{2},N-1}(s)
\bigtriangledown x_{1}(s)=0,
\qquad
0\leq k\leq n_{j}-2,
\qquad
j=1,\ldots,r.
\end{gather*}

From Lemma~\ref{LIK2}, we conclude that $\Delta K_{q,\vec{n}}^{\vec{p},\vec{\beta},N }(s) \in \mathbb{V}_{r}$. Hence, $\Delta K_{q,\vec{n}}^{\vec{p},\vec{\beta},N }(s)$ can be expressed as a linear combination of polynomials
$\left\{K_{q,\vec{n}-\vec{e}_{i}}^{\vec{p},\vec{\beta}_{i,1/q^{2}},N-1}(s) \right\}_{i=1}^{r}$, i.e.~
\begin{gather}
\Delta K_{q,\vec{n}}^{\vec{p},\vec{\beta},N }(s) =\sum\limits_{i=1}^{r}\xi_{i}K_{q,\vec{n}-\vec{e}_{i}}^{\vec{p},\vec{\beta}_{i,1/q^{2}},N-1}(s),
\qquad
\sum\limits_{i=1}^{r}\abs{\xi_{i}}>0.
\label{eq-delK2}
\end{gather}
Multiplying both sides of the equation~\eqref{eq-delK2} by 
$[s]_{q}^{(n_{k}-1)}\upsilon_{q}^{p_{k},\beta_k/q^{2},N-1}(s)
\bigtriangledown x_{1}(s)$ and using relations~\eqref{Orthog_CondK2} one has
\begin{multline}
\sum\limits_{s=0}^{N}\Delta K_{q,\vec{n}}^{\vec{p},\vec{\beta},N }(s)[s]_{q}^{(n_{k}-1)}\upsilon_{q}^{p_{k},\beta_k/q^{2},N-1}(s)\bigtriangledown x_{1}(s)\\
=\sum\limits_{i=1}^{r}\xi_{i}\sum\limits_{s=0}^{N}K_{q,\vec{n}-\vec{e}_{i}}^{\vec{p},\vec{\beta}_{i,1/q^{2}},N-1}(s)
[s]_{q}^{(n_{k}-1)}
\upsilon_{q}^{p_{k},\beta_k/q^{2},N-1}(s) \bigtriangledown x_{1}(s)
\\
=\xi_{k}\sum\limits_{s=0}^{N}K_{q,\vec{n}-\vec{e}_{k}}^{\vec{p},\vec{\beta}_{k,1/q^{2}},N-1}(s)[s]_{q}^{(n_{k}-1)}
\upsilon_{q}^{p_{k},\beta_k/q^{2},N-1}(s)\bigtriangledown x_{1}(s). \label{Ident_IK2}
\end{multline}
If we replace $[s]_{q}^{(k)}$ by $[s]_{q}^{(n_k-1)}$ in the left-hand side of equation~\eqref{inte-K12}, then equation~\eqref{Ident_IK2} transforms into 
\begin{multline}
\sum\limits_{s=0}^{N}\Delta K_{q,\vec{n}}^{\vec{p},\vec{\beta},N }(s)[s]_{q}^{(n_{k}-1)}\upsilon_{q}^{p_{k},\beta_k/q^{2},N-1}(s) \bigtriangledown x_{1}(s)
=-\sum\limits_{s=0}^{N}K_{q,\vec{n}}^{\vec{p},\vec{\beta},N }(s)\varphi_{k,n_{k}-1}(s)\upsilon_{q}^{p_{k},\beta_{k},N}(s) \bigtriangledown x_{1}(s)
\\
=\frac{q^{-1/2}\left(p_k\left(q^{n_{k} }-1\right)+1\right)}{q^{2(N-1)}(1-p_{k})[N]_q}\sum\limits_{s=0}^{\infty}K_{q,\vec{n}}^{\vec{p},\vec{\beta},N }(s)[s]_{q}^{(n_{k})}
\upsilon_{q}^{p_{k},\beta_{k},N}(s)\bigtriangledown x_{1}(s),
\label{eqcha-Kra2}
\end{multline}
in which we have used for the above transformation that $x(s)[s-1]_{q}^{(n_{k}-1)}=[s]_{q}^{(n_{k})}$. In addition,
\begin{equation*}
\varphi_{k,n_{k}-1}(s) 
=- \frac{q^{-1/2}\left(p_k\left(q^{n_{k}}-1\right)+1\right)}{q^{2(N-1)}p_{k}(1-p_k)[N]_q}[s]_{q}^{(n_{k})}+ \text{lower degree terms}.
\end{equation*}
From~\eqref{ROpqMKrav} one has
\begin{gather}
\frac{\left(p_k\left(q^{|\vec{n}|}-1\right)+1\right)}{q^{2(N-1)}p_{k}(1-p_k)[N]_q}\upsilon_{q}^{p_{k},\beta_{k},N}(s) K_{q,\vec{n}}^{\vec{p},\vec{\beta},N}(s) =-q^{|\vec{n}|-1/2}\nabla
\left(\upsilon_{q}^{p_{k},\beta_k/q^{2},N-1}(s) K_{q,\vec{n}-\vec{e}_{k}}^{\vec{p},\vec{\beta}_{k,1/q^{2}},N-1}(s)\right).
\label{eqcha-Kra12}
\end{gather}
Then, by substituting~\eqref{eqcha-Kra12} in the right-hand side of equation~\eqref{eqcha-Kra2} and using summation by parts, we get
\begin{multline*}
\sum\limits_{s=0}^{N}\Delta K_{q,\vec{n}}^{\vec{p},\vec{\beta},N}(s)
[s]_{q}^{(n_{k}-1)}\upsilon_{q}^{p_{k},\beta_k/q^{2},N-1}(s) \bigtriangledown x_{1}(s)\\
=-q^{|\vec{n}|-1}\frac{ p_{k}\left( q^{n_{k} }-1\right)+1}{ p_{k}\left( q^{|\vec{n}| }-1\right)+1}\sum\limits_{s=0}^{N}[s]_{q}^{(n_{k})}\nabla \big[\upsilon_{q}^{p_{k},\beta_k/q^{2},N-1}(s)K_{q,\vec{n}-\vec{e}_{k}}^{\vec{p},\vec{\beta}_{k,1/q^{2}},N-1}(s)\big]\bigtriangledown x_{1}(s)
\\
=q^{|\vec{n}|-1}\frac{ p_{k}\left( q^{n_{k} }-1\right)+1}{ p_{k}\left( q^{|\vec{n}| }-1\right)+1}\sum\limits_{s=0}^{N}K_{q,\vec{n}-\vec{e}_{k}}^{\vec{p},\vec{\beta}_{k,1/q^{2}},N-1}(s)\left(\Delta [s]_{q}^{(n_{k})}\right)\upsilon_{q}^{p_{k},\beta_k/q^{2},N-1}(s) \bigtriangledown x_{1}(s).
\end{multline*}
Since $\Delta [s]_{q}^{(n_{k})}=q^{3/2-n_{k}}[n_{k}]_{q}^{(1)}[s]_{q}^{(n_{k}-1)}$ we
finally obtain
\begin{multline*}
\sum\limits_{s=0}^{N}\Delta K_{q,\vec{n}}^{\vec{p},\vec{\beta},N}(s)
[s]_{q}^{(n_{k}-1)}\upsilon_{q}^{p_{k},\beta_k/q^{2},N-1}(s) \bigtriangledown x_{1}(s)
\\
=q^{|\vec{n}|-n_{k}+1/2}\left(\frac{ p_{k}\left( q^{n_{k} }-1\right)+1}{ p_{k}\left( q^{|\vec{n}| }-1\right)+1}\right)[n_{k}]_{q}^{(1)}\sum\limits_{s=0}^{N}K_{q,\vec{n}-\vec{e}_{k}}^{\vec{p},\vec{\beta}_{k,1/q^{2}},N-1}(s)
[s]_{q}^{(n_{k}-1)}\upsilon_{q}^{p_{k},\beta_k/q^{2},N-1}(s) \bigtriangledown x_{1}(s).
\end{multline*}
Comparing this equation with~\eqref{Ident_IK2} we obtain the coefficients in the expansion~\eqref{eq-delK2}
\begin{gather*}
\xi_{k}=q^{|\vec{n}|-n_{k}+1/2}\left(\frac{ p_{k}\left( q^{n_{k} }-1\right)+1}{ p_{k}\left( q^{|\vec{n}| }-1\right)+1}\right)[n_{k}]_{q}^{(1)}.
\end{gather*}
Therefore, the equation \refe{Rela_qK2} holds. 
\end{proof}

\begin{Theorem}\label{eqdifKra}
The $q$-Kravchuk multiple orthogonal polynomials $K_{q,\vec{n}}^{\vec{p},\vec{\beta},N}(s)$ satisfy the following 
$(r+1)$-order difference equation
\begin{multline}
\prod\limits_{i=1}^{r}\mathcal{D}_{q}^{p_{i},\beta_i/q^{2},N-1}\Delta K_{q,\vec{n}}^{\vec{p},\vec{\beta},N}(s)
\\
=-\sum\limits_{i=1}^{r}q^{\vert \vec{n}\vert -n_{i}+1}\left(\frac{ p_{i}\left( q^{n_{i} }-1\right)+1}{ p_{i}\left( q^{|\vec{n}| }-1\right)+1}\right)[n_{i}]_{q}^{(1)}
\prod\limits_{\substack{j=1
\\
j\neq i}}^{r}\mathcal{D}_{q}^{p_{j},\beta_j/q^{2},N-1}K_{q,\vec{n}}^{\vec{p},\vec{\beta},N}(s).
\label{q-DEquationK2}
\end{multline}
\end{Theorem}

\begin{proof}
Since the operators \eqref{D_raising_Krav} commute, we have
\begin{gather}
\prod\limits_{i=1}^{r}\mathcal{D}_{q}^{p_{i},\beta_i/q^{2},N-1}
=\left(\prod\limits_{\substack{j=1\\j\neq i}}^{r} \mathcal{D}_{q}^{p_{j},\beta_j/q^{2},N-1}\right) \mathcal{D}_{q}^{p_{i},\beta_i/q^{2},N-1}.
\label{intermK2}
\end{gather}
Using Lemma \ref{lemma_lower} and the raising operators \eqref{ROpqMKrav} in accordance with the above defined product of operators \eqref{intermK2} in the equation \eqref{Rela_qK2}, we obtain \eqref{q-DEquationK2}, i.e. 
\begin{multline*}
\prod\limits_{i=1}^{r}\mathcal{D}_{q}^{p_{i},\beta_i/q^{2},N-1}\Delta K_{q,\vec{n}}^{\vec{p},\vec{\beta},N}(s)\\
=\sum\limits_{i=1}^{r}q^{|\vec{n}|-n_{i}+1/2}\left(\frac{ p_{i}\left( q^{n_{i} }-1\right)+1}{ p_{i}\left( q^{|\vec{n}| }-1\right)+1}\right)[n_{i}]_{q}^{(1)}\prod\limits_{\substack{j=1\\j\neq i}}^{r}\mathcal{D}_{q}^{p_{j},\beta_j/q^{2},N-1}\left(\mathcal{D}_{q}^{p_{i},\beta_i/q^{2},N-1}K_{q,\vec{n}-\vec{e}_{i}}^{\vec{p},\vec{\beta}_{i,1/q^{2}},N-1}(s)\right)
\\
=-\sum\limits_{i=1}^{r}q^{\vert \vec{n}\vert -n_{i}+1}\left(\frac{ p_{i}\left( q^{n_{i} }-1\right)+1}{ p_{i}\left( q^{|\vec{n}| }-1\right)+1}\right)[n_{i}]_{q}^{(1)}\prod\limits_{\substack{j=1
\\
j\neq i}}^{r} \mathcal{D}_{q}^{p_{j},\beta_j/q^{2},N-1}K_{q,\vec{n}}^{\vec{p},\vec{\beta},N}(s).
\end{multline*}
\end{proof}

\begin{Remark}
When $r=1$ the above equation \refe{q-DEquationK2} coincides with the very important hypergeometric-type equation \refe{eqdif}, for the choice 
\begin{align}a_{2}(s)&=(q-1)x^{2}(s)+x(s),\notag \\
a_{1}(s) &=   -\dfrac{q^{\half}}{1-p}\left[p(q-1)+1\right]x(s)+
\frac{q^{\half}pq(q^{N}-1)}{1-p},\label{a2_a1_lambda}\\
\lambda_n&=q^{-\frac{n}{2} +1} [n]_q \dfrac{\left[p(q^{n}-1)+1\right]}{(1-p)}.\notag
\end{align}
\end{Remark}

Notice that for a function $f(s)$ defined on the discrete variable $s$, a straightforward calculation yields 
\begin{multline}
q^{-\vert \vec{n}\vert -1/2}\mathcal{D}_{q}^{p_{i}, \beta_{i}, N}f(s) =\frac{q^{-1}}{p_i\left( q^{|\vec{n}|-1}-1\right)+1}\bigg(p_{i}\left(x(N+1)-x(s)\right)-q\beta_{i}x(s) \bigg) f(s)\\
+\frac{\beta_{i}}{ p_{i} \left( q^{|\vec{n}|-1}-1\right)+1}x(s)
\bigtriangledown f(s).
\label{raising_different}
\end{multline}
Let $|\vec{n}|=n$, $p_i=p$, $\beta_i=(1-p)$ in \refe{a2_a1_lambda}, and combine the lowering operator \refe{lower} and \refe{raising_different} to obtain \refe{eqdif}.
\section{Connection with multiple Kravchuk polynomials}\label{limit-relations-Krav}

The limiting process when $q$ approaches $1$
transforms the $q$-algebraic relations studied in the present paper for multiple $q$-Kravchuk polynomials into the corresponding relations for multiple Kravchuk polynomials \cite{arvesu_vanAssche}.

\begin{Proposition}\label{limitMeixnerIProp}
The following relation hold:
\begin{equation}
\lim\limits_{q \rightarrow 1} K_{q,\vec{n}}^{\vec{p},\vec{\beta},N}(s) =K_{\vec{n}}^{\vec{p},N}(s),
\label{limitKravchuks}
\end{equation}
\end{Proposition}
where $K_{\vec{n}}^{\vec{p},N}(s)$ denotes the monic multiple Kravchuk polynomials \refe{RFK}.

\begin{proof} Let us transform the expression \refe{RFormulaKrav}, i.e.
\begin{gather*}
K_{q,\vec{n}}^{\vec{p},\vec{\beta},N}(s)
=\mathcal{G}_{q}^{\vec{n},\vec{p},N}\frac{\Gamma_{q}(N-s+1)\Gamma_{q}(s+1)}{q^{\binom{s}{2}}[N]_q!}
\prod\limits_{i=1}^{r}\left(\frac{\beta_i}{p_{i}}\right)^{s} \nabla^{n_{i}}\left(\frac{p_{i}q^{2n_i}}{\beta_i}\right)^{s}
\frac{q^{\binom{s}{2}}[N-|\vec{n}|]_q!}{\Gamma_{q}(N-\vert \vec{n} \vert-s+1)\Gamma_{q}(s+1)},
\end{gather*}
by using a $q$-analogue of \refe{Leibniz2} on the lattice $x(s)$ see \cite[formula~(3.2.29)]{nsu}
\begin{gather}
\nabla^{m}f(s)=q^{\binom{m+1}{2}/2-ms} \sum\limits_{k=0}^{m}{m\brack k}(-1)^k q^{\binom{m-k}{2}}f(s-k),
\label{qLeibniz2}
\end{gather}
where
\begin{gather*}
{m\brack k}=\frac{(q;q)_{m}}{(q;q)_{k}(q;q)_{m-k}},
\qquad
(a;q)_k=\prod\limits_{j=0}^{k-1}(1-aq^{j})
\quad k>0,
\quad
(a;q)_0=1.
\end{gather*}
Recall that $(a;q)_k$ denotes the $q$-analogue of the Pochhammer symbol in \refe{qLeibniz2} \cite{Gasper, Koek, nsu}. 

From \refe{RFormulaKrav} and \refe{qLeibniz2} we have
\begin{multline}
K_{q,\vec{n}}^{\vec{p},\vec{\beta},N}(s)
=\mathcal{G}_{q}^{\vec{n},\vec{p},N}q^{\sum_{j=1}^{r}\binom{n_{j}+1}{2}/2}
\dfrac{\Gamma_q (N-s+1)\Gamma_q (s+1)}{q^{\binom{s}{2}}[N]_q!}
\\
\times\sum\limits_{k_r=0}^{n_r}\cdots\sum\limits_{k_1=0}^{n_1}(-1)^{|\vec{k}|} {n_r\brack k_r} \cdots {n_1\brack k_1}
\frac{q^{|\vec{n}|s+\sum_{j=1}^{r}\binom{n_j-k_j}{2}-2k_j n_{j}-\sum_{j=2}^{r}k_{j}|\vec{n}|_{j}}}
{(p_r/\beta_r)^{k_{r}}(p_{r-1}/\beta_{r-1})^{k_{r-1}}\cdots(p_1/\beta_1)^{k_{1}}}
\\
\times 
\frac{ q^{\binom{s-|\vec{k}|}{2}}
[N-|\vec{n}|]_q!}{\Gamma_q(N-|\vec{n}|-s+|\vec{k}|+1)\Gamma_q(s-|\vec{k}|+1)}.
\label{multi_krav_4_limit}
\end{multline}
Calculating the limit as $q$ approaches $1$ in \refe{multi_krav_4_limit} we have
\begin{multline}
\lim\limits_{q \rightarrow 1}K_{q,\vec{n}}^{\vec{p},\vec{\beta},N}(s)
=\left( -N\right) _{\left\vert \vec{n}\right\vert }\left( \prod\limits_{i=1}^{r}p_{i}^{n_{i}}\right)
\dfrac{\Gamma (N-s+1)\Gamma (s+1)}{[N]!}
\\
\times\sum\limits_{k_r=0}^{n_r}\cdots\sum\limits_{k_1=0}^{n_1}(-1)^{|\vec{k}|} 
\binom{n_r}{ k_r} \cdots \binom{n_1}{ k_1}
\frac{1}
{(p_r/\beta_r)^{k_{r}}(p_{r-1}/\beta_{r-1})^{k_{r-1}}\cdots(p_1/\beta_1)^{k_{1}}}
\\
\times 
\frac{[N-|\vec{n}|]!}{\Gamma(N-|\vec{n}|-s+|\vec{k}|+1)\Gamma(s-|\vec{k}|+1)}.
\label{multi_krav_limit}
\end{multline}
Now, using \refe{Leibniz2}, \refe{operator_Krav_Rodrig} and noticing that $\beta_{i}=1-p_{i}$ we transform the right-hand side of expression \refe{multi_krav_limit} into the equation
\begin{align*}
\lim\limits_{q \rightarrow 1}K_{q,\vec{n}}^{\vec{p},\vec{\beta},N}(s)
&=\left( -N\right) _{\left\vert \vec{n}\right\vert }\left( \prod\limits_{i=1}^{r}p_{i}^{n_{i}}\right) \frac{\Gamma(s+1)\Gamma (N-s+1)}{N!}\\
&\times\prod_{i=1}^{r}\left( \frac{1-p_{i}}{p_{i}}\right)^{s}\bigtriangledown ^{n_{i}}\left( \frac{p_{i}}{1-p_{i}}\right)^{s}
\frac{\left( N-\left\vert \vec{n}\right\vert \right) !}{\Gamma (s+1)\Gamma (N-\left\vert \vec{n}\right\vert-s+1)}\\
&=K_{\vec{n}}^{\vec{p},N}(s),
\end{align*}
which proves \refe{limitKravchuks}.

\end{proof}
\begin{Remark}
The difference equation for multiple $q$-Kravchuk polynomials \refe{q-DEquationK2} transforms into \refe{opdi-1Kra} as $q$ approaches to $1$. Indeed, we follow the above description: Use \refe{qLeibniz2}, the raising operators \refe{D_raising_Krav}, and compute the limit $q\to1$ to suitably transform equation \refe{q-DEquationK2}. Then, conveniently rewriting the obtained expression in the uniform lattice $s=0,1,\dots,N$ and with the help of expression \refe{Leibniz2} and the raising operators \refe{raisingxKra} we get the difference equation \refe{opdi-1Kra} for multiple Kravchuk polynomials.
\end{Remark}

Finally, we summarize the connection between the following four polynomial families involved in the present paper (Kravchuk polynomials of different types) 

\begin{equation*}
  \begin{CD}
    K_{q,\vec{n}}^{\vec{p},\vec{\beta},N}(s)  @>q\rightarrow 1>>
    K_{\vec{n}}^{\vec{p},N}(s) \\
    @Vr=1VV @Vr=1VV \\
    K_{q,n}^{p,N}(s) @>q\rightarrow 1>> K_{n}^{p,N}(s) 
  \end{CD}
\end{equation*}

Furthermore,

\begin{equation*}
  \begin{CD}
    \fbox{Diff. Eq. \refe{q-DEquationK2} on $x(s)$ for $K_{q,\vec{n}}^{\vec{p},\vec{\beta},N}(s)$}  @>q\rightarrow 1>>
    \fbox{Diff. Eq. \refe{opdi-1Kra} on $s$ for $K_{\vec{n}}^{\vec{p},N}(s)$} \\
    @Vr=1VV @Vr=1VV \\
    \fbox{Hypergeometric-type Eq. \refe{eqdif} for $K_{q,n}^{p,N}(s)$} @>q\rightarrow 1>> \fbox{Hypergeometric-type Eq. \refe{eqdif} for $K_{n}^{p,N}(s)$} 
  \end{CD}
\end{equation*}


\begin{thebibliography}{999}

\bibitem{alv_arv}
{\'A}lvarez-Nodarse R., Arves{\'u} J., On the {$q$}-polynomials in the
exponential lattice {$x(s)=c_1q^s+c_3$}, \href{http://dx.doi.org/10.1080/10652469908819236}{\textit{Integral Transform. Spec.
Funct.}} \textbf{8} (1999), 299--324.

\bibitem{angelesco}
Angelesco A., Sur l'approximation simultan\'ee de plusieurs int\'egrales d\'efinies, C.R. Acad. Sci, Paris \textbf{167} (1918), 629--631.

\bibitem{Aptekarev1} Aptekarev A.I., Multiple orthogonal polynomials, {\textit{J. Comput. Appl. Math.}} \textbf{99} (1998), 423--447.

\bibitem{arvesu-qalgebras} 
Arves\'u J. Quantum algebras $su_{q}(2)$ and $su_{q}(1,1)$ associated with certain $q$-Hahn polynomials: A revisited approach, Electron. Trans. Numer. Anal., 2006; 24: 24--44.

\bibitem{arvesu-qHahn}
Arves\'u J., On some properties of $q$-Hahn multiple orthogonal
polynomials, \href{http://dx.doi.org/10.1016/j.cam.2009.02.062}{\textit{J.~Comput. Appl. Math.}} \textbf{233} (2010), 1462--1469.

\bibitem{arvesu_vanAssche}
Arves{\'u} J., Coussement J., Van~Assche W., Some discrete multiple orthogonal
polynomials, \href{http://dx.doi.org/10.1016/S0377-0427(02)00597-6}{\textit{J.~Comput. Appl. Math.}} \textbf{153} (2003), 19--45.

\bibitem{arvesu-esposito}
Arves{\'u} J., Esposito C., A high-order $q$-difference equation for
{$q$}-{H}ahn multiple orthogonal polynomials, \href{http://dx.doi.org/10.1080/10236198.2010.524211}{\textit{J.~Difference Equ.
Appl.}} \textbf{18} (2012), 833--847.


\bibitem{arvesu-andys}
Arves\'u J., Ram\'irez-Aberasturis A.M., On the $q$-Charlier multiple orthogonal polynomials, {\textit{SIGMA}} \textbf{11} (2015), 026.

\bibitem{arvesu-andys1}
Arves\'u J., Ram\'irez-Aberasturis A.M., On some algebraic properties for $q$-Meixner multiple orthogonal polynomials of the first kind, {\textit{arXiv:1604.04920}} (2016), 1--10.

\bibitem{arvesu-andys2}
Arves\'u J., Ram\'irez-Aberasturis A.M., Multiple Meixner polynomials on a non-uniform lattice,
\textit{Mathematics}, \textbf{8}, (2020), 1460; doi:10.3390/math8091460.

\bibitem{Coussement-Assche}
Coussement J., Van Assche W., Differential equations for multiple orthogonal polynomials with respect to classical weights, {\textit{J. Phys. A: Math. Gen.}} \textbf{39} (2006), 3311--3318.

\bibitem{cse} C. Campigotto, Yu. F. Smirnov, and S. G. Enikeev: 
$q$-Analogue of the Kravchuk and Meixner orthogonal polynomials, 
\textit{J. of Comput. and Appl. Math.},  {\bf 57}, (1995), 87--97.

\bibitem{floreanini-letourneux-vinet}
Floreanini R., LeTourneux J., Vinet L., Quantum mechanics and polynomials of a discrete variable, \href{https://doi.org/10.1006/aphy.1993.1072}{\textit{Ann. Phys.}} \textbf{226} (1993), 331--349.

\bibitem{Gasper}
Gasper G., Rahman M., Basic hypergeometric series, \href{http://dx.doi.org/10.1017/CBO9780511526251}{\textit{Encyclopedia of
Mathematics and its Applications}}, Vol.~96, 2nd ed., Cambridge University
Press, Cambridge, (2004).

\bibitem{Gonchar}
Gonchar A.A., Rakhmanov E.A., Sorokin V.N., Hermite-Pad\'e approximants for systems of Markov-type functions, {\textit{Mat. Sb.}} \textbf{188} (1997), 33--58.

\bibitem{Hane}
Haneczok M., Van~Assche W., Interlacing properties of zeros of multiple
orthogonal polynomials, \href{http://dx.doi.org/10.1016/j.jmaa.2011.11.077}{\textit{J.~Math. Anal. Appl.}} \textbf{389} (2012),

\bibitem{Kalyagin1}
Kalyagin V.A., On a class of polynomials defined by two orthogonality relations, \href{http://mi.mathnet.ru/eng/msb2512}{\textit{Mat. Sb.}} \textbf{110} (1979), 609--627.

\bibitem{Kalyagin}
Kalyagin(Kaliaguine) V.A., Higher order difference operator's spectra characteristics
and the convergence of the joint rational approximations, {\textit{Dokl. Akad. Nauk}} \textbf{340} (1995), 15--17.

\bibitem{Kalyagin2}
Kaliaguine(Kalyagin) V.A., Ronveaux A., On a system of classical polynomials of simultaneous orthogonality, \href{https://doi.org/10.1016/0377-0427(94)00129-4}{\textit{J. Comput. Appl. Math.}} \textbf{67} (1996), 207--217.

\bibitem{Koek}
Koekoek R., Lesky P.A., Swarttouw R. F., Hypergeometric orthogonal polynomials
and their $q$-analogues, \href{http://dx.doi.org/10.1007/978-3-642-05014-5}{\textit{Springer Monographs in Mathematics}},
Springer-Verlag, Berlin, (2010).

\bibitem{kravchuk1} M. Kravchuk, On interpolation by means of
orthogonal polynomials,  \textit{Zap. Kiev. Sil'sk.-Gospod. Inst.}, 
{\bf 4}, (1929), 21--28.

\bibitem{kravchuk2} M. Kravchuk, On orthogonal polynomials
associated with trials of inverted and noninverted circles, \textit{ 
Zap. Fiz. Mat. Viddilu Ukr. Akad. Nauk.}, {\bf 5}, (1931), 19--48.

\bibitem{lee}
Lee D.W., Difference equations for discrete classical multiple orthogonal polynomials, \href{http://dx.doi.org/10.1016/j.jat.2007.06.002}{\textit{J.~Approx. Theory}} \textbf{150} (2008), 132--152.

\bibitem{Mahler}
Mahler K., Perfect systems, \href{http://eudml.org/doc/88959}{\textit{Compos. Math.}} \textbf{19} (1968), 95--166.

\bibitem{miki-tsujimoto-vinet-zhedanov}
Miki H., Tsujimoto S., Vinet L., Zhedanov A., An algebraic model for the multiple
Meixner polynomials of the first kind, {\textit{J. Phys.~A: Math. Theor. }} \textbf{45} (2012), 325205 (11 pp).

\bibitem{miki-vinet-zhedanov}
Miki H., Vinet L., Zhedanov A., Non-{H}ermitian oscillator {H}amiltonians and
multiple {C}harlier polynomials, \href{http://dx.doi.org/10.1016/j.physleta.2011.10.038}{\textit{Phys. Lett.~A}} \textbf{376} (2011),
65--69.

\bibitem{nda_van_assche} Ndayiragije F., Van Assche W., Multiple Meixner polynomials and non-Hermitian oscillator Hamiltonians, \textit{J. Phys. A} \textbf{46} (2013), 505201 (17 pp).

\bibitem{nsu}
Nikiforov A.F., Suslov S.K., Uvarov V.B., Classical orthogonal polynomials of a
discrete variable, \href{http://dx.doi.org/10.1007/978-3-642-74748-9}{\textit{Springer Series in Computational Physics}}, Springer-Verlag,
Berlin, (1991).

\bibitem{kn_Nikishin} 
Nikishin E.M., On simultaneous Pad\'{e} approximations, {\textit{Mat. Sb.}} \textbf{113} (1980), 499--519; English transl. Math. USSR Sb., \textbf{41} (1982).

\bibitem{Nikishin}
Nikishin E.M., Sorokin V.N., Rational approximations and orthogonality,
\textit{Translations of Mathematical Monographs}, Vol.~92, Amer. Math. Soc.,
Providence, RI, (1991).

\bibitem{Sorokin}
Sorokin V.N., A generalization of classical orthogonal polynomials and the convergence of simultaneous Pad\'e approximants, {\textit{Trudy Sem. Im. I. G. Petrovsk.}} \textbf{11} (1986), 125--165.

\bibitem{Sorokin3}
Sorokin V.N., Simultaneous Pad\'e approximants for finite and infinite intervals, \href{http://mi.mathnet.ru/eng/ivm8510}{\textit{Izv. Vyssh. Uchebn. Zaved. Mat.}} \textbf{8} (1984), 45--52.

\bibitem{Assche_diff-eq}
Van~Assche W., Difference equations for multiple {C}harlier and {M}eixner
polynomials, in Proceedings of the {S}ixth {I}nternational {C}onference on
{D}ifference {E}quations, CRC, Boca Raton, FL, (2004), 549--557.

\end{thebibliography}
\end{document}